\documentclass[11pt]{article}
\usepackage{latexsym}
\usepackage{amsthm,amsmath,amssymb,amstext}

\usepackage{graphicx} 
\usepackage{psfrag,color}
\usepackage{slashbox,enumerate}

\numberwithin{equation}{section}

\def\XXint#1#2#3{{\setbox0=\hbox{$#1{#2#3}{\int}$}
     \vcenter{\hbox{$#2#3$}}\kern-.5\wd0}}

\newcommand{\adj}[1]{p_h(y_h(#1))}

\newtheorem{Theorem}{Theorem}[section]

\newtheorem{Lemma}[Theorem]{Lemma}

\theoremstyle{definition}

\newtheorem{Algorithm}[Theorem]{Algorithm}

\newtheorem{Remark}[Theorem]{Remark}

\newtheorem{Example}[Theorem]{Example}

\newtheorem*{Examples}{Examples}

\newtheorem{Assumption}[Theorem]{Assumption}

\def\mathref#1{\ifmmode\mathrm{(\ref{#1})}\else(\ref{#1})\fi}

\def\rhx2{\sqrt{1+ r_{h,x}^2}}

\title{Variational discretization and semi-smooth Newton methods; implementation,
  convergence and globalization in pde constrained optimization with control constraints}
\author{Michael Hinze\footnote{Schwerpunkt Optimierung und Approximation,
Universit\"at Hamburg, Bundesstra{\ss}e 55, 20146 Hamburg, Germany.} \& Morten
Vierling\footnote{Schwerpunkt Optimierung und Approximation,
Universit\"at Hamburg, Bundesstra{\ss}e 55, 20146 Hamburg, Germany.}}

\date{\today}

\begin{document}

\maketitle

\begin{center}
{\bf \LARGE  }
\end{center}

{\small {\bf Abstract:} When combining the numerical concept of variational discretization introduced
  in \cite{H03,H05} and semi-smooth Newton methods for the numerical solution
  of pde constrained optimization with control constraints \cite{HIK03,UM03} special emphasis has to be
  taken on the implementation, convergence and globalization of the numerical
  algorithm. In the present work we address all these issues. In particular we
prove fast local convergence of the algorithm and propose two different globalization
strategies which are applicable in many practically relevant mathematical
settings. We illustrate our analytical and algorithmical findings by numerical experiments.} \\[2mm]

\pagenumbering{arabic}

\section{Introduction and mathematical setting}\label{sec:intro}
We are interested in the numerical treatment of the following control problem
\begin{equation}\label{H:MP}
(\mathbb{P})\quad \left\{\begin{array}{l}
\min_{(y,u)\in Y\times U_{ad}} J(y,u) := \frac{1}{2}\|y-z\|_{L^2(\Omega)}^2 + \frac{\alpha}{2} \|u\|_U^2 \\
\mbox{s.t.}\\
\begin{array}{rcll} 
-\Delta y &= &Bu & \mbox{ in } \Omega,\\
y & = & 0 & \mbox{ on } \partial\Omega\,.
\end{array}
\end{array}\right.
\end{equation}
Here, $\Omega \subset \mathbb{R}^d\, (d\ge 1)$ denotes an open, bounded sufficiently
smooth (polyhedral) domain. Given some Hilbert space $U$ and some closed, convex admissible set $U_{ad}\subset U$ for the controls and a linear, continuous control operator $B: U\rightarrow H^{-1}(\Omega)$, the states lie in $Y := H^1_0(\Omega)$. Let us note that also additional state constraints could be included into our problem setting, as done in \cite{DeHi06} and \cite{DH07}, and also more general (linear) elliptic or parabolic state equations. However, all structural issues discussed in the present work are induced by the control constraints, hence to keep the exposition as simple as possible state constraints are not considered here.

Typical configurations of $\mathbb{P}$ are
\begin{Examples}$\mbox{}$
\begin{enumerate}[(i)]
\item $U := \mathbb{R}^m, \,Y=H_0^1(\Omega),\,B: \mathbb{R}^m \rightarrow H^{-1}(\Omega), \ Bu := \sum\limits_{j=1}^{m} u_j
 F_j, \ F_j \in H^{-1}(\Omega),\, U_{\text{ad}} := \{v \in \mathbb{R}^m; a_j \le v_j \le b_j\},\,a,b\in\mathbb{R}^m,\, a < b.$
 
\item  $U := L^2(\Omega), \,Y=H_0^1(\Omega),\,B=\imath : L^2(\Omega) \rightarrow H^{-1}(\Omega)$, $\imath$ being the canonical injection, $ U_{\text{ad}} := \{v \in L^2(\Omega); a\le v\le b\},\,a,b\in L^\infty(\Omega),\, a < b.$
 
 \end{enumerate}
\end{Examples}
\begin{Remark}
One may as well consider elliptic equations with Neumann boundary control, 
\[
\begin{array}{rcll} 
-\Delta y +y&= &0 & \mbox{ in } \Omega,\\
\partial_\eta y & = & u & \mbox{ on } \partial\Omega\,,
\end{array}
\]
thus setting $U := L^2(\Gamma), \,Y=H^1(\Omega)$, 
$ U_{\text{ad}} := \{v \in L^2(\Gamma); a\le v\le b\},\,a,b\in L^\infty(\Gamma),\, a < b.$
\end{Remark}
Problem $\mathbb{P}$ admits a unique solution $(y,u) \in Y\times
U_{\text{ad}}$, and can equivalently be rewritten as the optimization problem
\begin{equation}\label{H:OP}
\min_{u \in U_{\text{ad}}}{\hat J(u)}
\end{equation}
for the reduced functional $\hat J(u) := J(y(u),u) \equiv J(SBu,u)$ over the set $U_{\text{ad}}$, where $S: Y^* \rightarrow Y$ denotes the (continuous) solution operator associated with $-\Delta$. We further know that the first order necessary (and here also sufficient) optimality conditions take the form
\begin{equation}\label{H:VI}
\langle\hat J'(u),v-u\rangle_{U^*,U} \ge 0 \mbox{ for all } v \in U_{\text{ad}}
\end{equation}
where $\hat J'(u) = \alpha u + B^*S^*(SBu-z) \equiv \alpha u + B^* p$, with $p := S^*(SBu-z)$ denoting the adjoint variable. The function $p$ in our setting satisfies
\begin{equation}
  \label{H:adjoint}
  \begin{array}{rcll} 
-\Delta p &= &y-z & \mbox{ in } \Omega, \\
p & = & 0 & \mbox{ on } \partial\Omega.
\end{array}
\end{equation}
For the numerical treatment of problem \eqref{H:MP} it is convenient to
rewrite \eqref{H:VI} for $\sigma > 0$ arbitrary in form of the following non--smooth operator equation;
\begin{equation*}
  u = P_{U_{\text{ad}}}\left(u-\sigma \nabla \hat J(u)\right) \stackrel{\sigma=1/\alpha}{\equiv} P_{U_{\text{ad}}}\left(-\frac{1}{\alpha}R^{-1}B^*p\right),
\end{equation*}
with the Riesz isomorphism $R: U\rightarrow U^*$ and the gradient $\nabla \hat J(u) = R^{-1}\hat J'(u)$.
\section{Finite element discretization}\label{H:sec:FEsemiD}
To discretize $(\mathbb{P})$ we concentrate on Finite Element approaches and make the following assumptions.
\begin{Assumption}\label{H:FEAssumptions}$\mbox{}$\\
$\Omega \subset \mathbb{R}^d$ denotes a polyhedral domain, $\bar \Omega = \cup_{j=1}^{nt} \bar T_j$ with admissible quasi-uniform sequences of partitions $\{T_j\}_{j=1}^{nt}$ of $\Omega$, i.e. with $h_{nt}:= \max_j\mbox{diam }T_j$ and $\sigma_{nt}:= \min_j\{\sup\mbox{diam }K; K \subseteq T_j\}$ there holds $c \le \frac{h_{nt}}{\sigma_{nt}} \le C$ uniformly in $nt$ with positive constants $0 < c \le C < \infty$ independent of $nt$. We abbreviate $T_h := \{T_j\}_{j=1}^{nt}$.
\end{Assumption}
For $k \in \mathbb{N}$ we set
\begin{multline*}
W_h := \{v \in C^0(\bar\Omega); v_{\mid_{T_{j}}} \in \mathbb{P}_k(T_j) \mbox{ for all } 1\le j \le nt\} =: \langle \phi_1,\dots,\phi_{ng}\rangle, \mbox{ and} \\ Y_h:= \{v \in W_h, v_{\mid_{\partial\Omega}} = 0\} =: \langle \phi_1,\dots,\phi_{n}\rangle  \subseteq Y,
\end{multline*}
with some $0 < n < ng$. The resulting Ansatz for $y_h$ then is of the form
$y_h = \sum\limits_{i=1}^{n} y_i \phi_i$. Now we approximate problem $(\mathbb{P})$ by
\begin{equation}\label{H:MPd}
(\mathbb{P}_h)\quad \left\{\begin{array}{l}
\min_{(y_h,u)\in Y_h\times U_{\text{ad}}} J(y_h,u) := \frac{1}{2}\|y_h-z\|_{L^2(\Omega)}^2 + \frac{\alpha}{2} \|u\|_U^2 \\
\mbox{s.t.}\\
\begin{array}{rcll} 
a(y_h,v_h) &= &\langle Bu,v_h\rangle_{Y^*,Y} & \mbox{ for all } v_h \in Y_h,
\end{array}\\
\end{array}\right.
\end{equation}
where $a(y,v) := \int_\Omega \nabla y \nabla v dx$ denotes the bilinear form associated with $-\Delta$. Problem
$(\mathbb{P}_h)$ admits a unique solution $(y_h,u) \in Y_h \times
U_{\text{ad}}$ and, as above, can equivalently be rewritten as the optimization problem
\begin{equation}\label{H:OPd}
\min_{u \in U_{\text{ad}}}{\hat J_h(u)}
\end{equation}
for the discrete reduced functional $\hat J_h(u) := J(y_h(u),u) \equiv
J(S_hBu,u)$ over the set $U_{\text{ad}}$, where $S_h: Y^* \rightarrow Y_h \subset Y $
denotes the solution operator associated with the finite element
discretization of $-\Delta$. The first order necessary (and here also sufficient) optimality conditions take the form
\begin{equation}\label{H:VId}
\langle\hat J_h'(u_h),v-u_h\rangle_{U^*,U} \ge 0 \mbox{ for all } v \in U_{\text{ad}}
\end{equation}
where $\hat J_h'(v) = \alpha v + B^*S_h^*(S_hBv-z) \equiv \alpha u + B^* p_h$, with $p_h := S_h^*(S_hBu-z)$ denoting the adjoint variable. The function $p_h$ in our setting satisfies
\begin{equation}
  \label{H:adjointd}
  a(v_h,p_h) = \langle y_h-z,v_h\rangle_{Y^*,Y} \mbox{ for all } v_h \in Y_h.
\end{equation}
Analogously to \eqref{H:VI}, for $\sigma > 0$ arbitrary, we have
\begin{equation}
  \label{H:OEd}
  u_h = P_{U_{\text{ad}}}\left(u_h-\sigma \nabla\hat J_h(u_h)\right) \stackrel{\sigma=1/\alpha}{\equiv} P_{U_{\text{ad}}}\left(-\frac{1}{\alpha}R^{-1}B^*p_h\right)\,.
\end{equation}
\begin{Remark}{\rm
Problem \eqref{H:MPd} is still infinite--dimensional in that the control space
is not discretized. This is reflected through the appearance of the projector
$ P_{U_{\text{ad}}}$ in \eqref{H:OEd}. The numerical challenge now consists in
designing numerical solution algorithms for problem \eqref{H:MPd} which are
implementable, and which reflect the infinite--dimensional structure of the
{\it discrete} problem \eqref{H:MPd} \cite{H03,H05}.
}
\end{Remark}

Next let us investigate the error $\|u-u_h\|_U$ between the solutions
$u$ of \eqref{H:OP} and $u_h$ of \eqref{H:OPd}, compare \cite{HinzePinnauUlbrich2009}.
\begin{Theorem}
Let $u$ denote the unique solution of \eqref{H:OP}, and $u_h$ the unique
solution of \eqref{H:OPd}. Then there holds
\begin{multline}\label{H:EE1}
\alpha\|u-u_h\|_U^2 + \frac{1}{2}\|y(u)-y_h\|^2 \le \langle B^*(p(u) -
\tilde p_h(u)),u_h-u\rangle_{U*,U} + \frac{1}{2}\|y(u)-y_h(u)\|_{L^2(\Omega)}^2,
\end{multline}
where $\tilde p_h(u) := S_h^*(SBu-z)$, $y_h(u) := S_hBu$, and $y(u) := SBu$.
\end{Theorem}
\noindent
\begin{proof} Since \eqref{H:OPd} is an optimization problem defined on
all of $U_{\text{ad}}$, the unique solution $u$ of \eqref{H:OP} is an
admissible test function in \eqref{H:VId}. Let us emphasize, that this
is different for approaches, where the control space is discretized
explictly. In this case we may only expect that $u_h$ is an admissible
test function for the continuous problem (if ever). So let us test
\eqref{H:VI} with $u_h$, and \eqref{H:VId} with $u$, and then add the
resulting variational inequalities. This leads to
\[
\left\langle \alpha(u-u_h) + B^*S^*(SBu-z) -
B^*S_h^*(S_hBu_h-z),u_h-u\right\rangle_{U^*,U} \ge 0.
\]
This inequality is equivalent to
\[
\alpha\|u-u_h\|_U^2 \le \left\langle B^*(p(u)-\tilde p_h(u)) +
B^*(\tilde p_h(u)-p_h(u_h)),u_h-u\right)\rangle_{U^*,U}.
\]
Let us investigate the second addend on the right hand side of this
inequality. By definition of the adjoint variables there holds
\begin{multline*}
\left\langle B^*(\tilde p_h(u)-p_h(u_h),u_h-u\right\rangle_{U^*,U} =
\langle \tilde p_h(u)-p_h(u_h),B(u_h-u)\rangle_{Y,Y^*} = \\
=a(y_h-y_h(u),\tilde p_h(u)-p_h(u_h)) =
\int\limits_{\Omega}(y_h(u_h)-y_h(u))(y(u)-y_h(u_h)) dx = \\ =-\|y_h-y\|_{L^2(\Omega)}^2
+ \int\limits_{\Omega}(y-y_h)(y-y_h(u)) dx \le -\frac{1}{2}\|y_h-y\|_{L^2(\Omega)}^2 +
\frac{1}{2}\|y-y_h(u)\|_{L^2(\Omega)}^2
\end{multline*}
so that the claim of the theorem follows.
\end{proof}
What can we learn from Theorem \ref{H:EE1}? It tells us that an error
estimate for $\|u-u_h\|_U$ is at hand, if
\begin{itemize}
 \item an error estimate for $\|R^{-1}B^*(p(u)-\tilde p_h(u)\|_U$ is
available, and
 \item an error estimate for $\|y(u)-y_h(u)\|_{L^2(\Omega)}$ is available.
\end{itemize}
\begin{Remark}
The error $\|u-u_h\|_U$ between the solution $u$ of problem \eqref{H:OP}
and $u_h$ of \eqref{H:OPd} is completely determined by the approximation
properties of the discrete solution operators $S_h$ and $S_h^*$.
\end{Remark}
%

\section{Semi-smooth Newton algorithm}


In the following we restrict our considerations to the practically relevant case of the second example given in Section \ref{sec:intro}, i.e. we set $U=L^2(\Omega)$, $Y=H^1_0(\Omega)$, $U_{ad} =\{v\in L^2(\Omega);\,a\le v \le b\}$ with $a,b\in L^\infty(\Omega)$, $b-a>\sigma>0$ and $\sigma\in \mathbb R$. Also the control operator is the injection $\imath:L^2(\Omega)\rightarrow Y^*$, hence the adjoint $B^*=\imath^*$ is the injection from $Y$ into $L^2(\Omega)$. Below, the operators $B$, $B^*$ and $R$ are omitted for notational convenience. The variationally discretized problem associated to $(\mathbb P)$ then reads
\begin{equation*}
(\mathbb P_h)\quad \left\{\begin{array}{l}
\min_{(y_h,u)\in Y\times L^2(\Omega)} J(y,u) := \frac{1}{2}\|y-z\|_{L^2(\Omega)}^2 + \frac{\alpha}{2} \|u\|_{L^2(\Omega)}^2 \\
\mbox{s.t.}\\
\begin{array}{rcll} 
a(y_h,v_h) &= &\langle u,v_h\rangle_{{L^2(\Omega)}} & \mbox{ for all } v_h \in Y_h
\end{array}\\
\mbox{and}\\
a\le u \le b \textrm{, a.e. in $\Omega$}\,.
\end{array}\right.
\end{equation*}%
To apply the semi-smooth Newton algorithm proposed in the following, the bounds are required to be elements of the finite element space $Y_h$. Let therefore $a_h,b_h\in Y_h$ be obtained from $a,b$ by interpolation or projection and let us consider the problem
\begin{equation*}
(\mathbb P_{hh})\quad \left\{\begin{array}{l}
\min_{(y_{hh},u)\in Y\times {L^2(\Omega)}} J(y,u) := \frac{1}{2}\|y-z\|_{L^2(\Omega)}^2 + \frac{\alpha}{2} \|u\|_{L^2(\Omega)}^2 \\
\mbox{s.t.}\\
\begin{array}{rcll} 
a(y_{hh},v_h) &= &\langle u,v_h\rangle_{L^2(\Omega)} & \mbox{ for all } v_h \in Y_h
\end{array}\\
\mbox{and}\\
a_h\le u \le b_h \textrm{, a.e. in $\Omega$},
\end{array}\right.
\end{equation*}
It is clear that for $h>0$ small enough the admissible set $a_h\le u\le b_h$ is non empty, if $a_h,b_h\stackrel{h\rightarrow 0}{\longrightarrow} a,b$ uniformly, say which can be guaranteed for sufficiently smooth bounds $a,b$ and $a_h=I_ha$, $b_h=I_hb$, with $I_h$ denoting the Lagrange interpolation operator or the $L^2$-projection. In this case problem $(\mathbb P_{hh})$ admits a unique solution $(u_{hh},y_{hh})$. Let us assume, that this solution exists.
\begin{Lemma}[Perturbed Bounds]\label{V:LemmaPertBnds}
The solutions $(y_{hh}, u_{hh})$ and $( y_{h}, u_{h})$ of $(\mathbb P_{hh})$ and $(\mathbb P_{h})$ satisfy the estimate
\[
\| u_{hh}-u_{h}\|_{L^2(\Omega)} \le\left(1+\frac{1}{\alpha}\|S_h\|^2\right)(\|a-a_h\|_{L^2(\Omega)}+\|b-b_h\|_{L^2(\Omega)})\,.
\]
\end{Lemma}
\begin{proof}
Let
\[
u_h^p(\omega) = P_{[a_h(\omega),b_h(\omega)]}\left(-\frac{1}{\alpha}[S^*_h(S_h u_h-z)](\omega)\right)\,.
\]
Then by (\ref{H:OEd}) there holds
\begin{equation}\label{V:BndsEst}
\|u_h^p-u_h\|_{L^2(\Omega)}\le \|a_h-a\|_{L^2(\Omega)}+\|b_h-b\|_{L^2(\Omega)}\,.
\end{equation}
Since $u_h^p$ is admissible for $\mathbb P_{hh}$ 
we have
\[
\langle u_{hh}+\frac{1}{\alpha}S^*_h(S_h u_{hh}-z),u_{h}^p- u_{hh}\rangle_{L^2(\Omega)} \ge 0
\]
while the definition of $u_h^p$ gives
\[
\langle u_{h}^p+\frac{1}{\alpha}S^*_h(S_h u_{h}-z), u_{hh}-u_{h}^p\rangle_{L^2(\Omega)} \ge 0
\]
since $ u_{hh}$ lies between $a_h$ and $b_h$. Adding these inequalities leads to
\[
\begin{split}
 \|u_{h}^p- u_{hh}\|^2_{L^2(\Omega)} &\le \frac{1}{\alpha} \langle S^*_hS_h( u_{h}- u_{hh}), u_{hh}-u_{h}^p\rangle_{L^2(\Omega)}\\
 & =\frac{1}{\alpha} \langle S^*_hS_h( u_{h}-u_{h}^p), u_{hh}-u_{h}^p\rangle_{L^2(\Omega)} +
 \frac{1}{\alpha} \langle S^*_hS_h(u_{h}^p- u_{hh}), u_{hh}-u_{h}^p\rangle_{L^2(\Omega)}
 \end{split}
\]
and finally we have
\[
\begin{split}
\|u_{h}^p- u_{hh}\|_{L^2(\Omega)}^2 + \frac{1}{\alpha}\|S_h(u_{h}^p- u_{hh})\|_{L^2(\Omega)}^2 &\le \frac{1}{\alpha}\langle S^*_hS_h( u_{h}-u_{h}^p), u_{hh}-u_{h}^p\rangle_{L^2(\Omega)}\\
&\le \frac{1}{\alpha}\|S_h\|^2\| u_{h}-u_{h}^p\|_{L^2(\Omega)}\| u_{hh}-u_{h}^p\|_{L^2(\Omega)}
\end{split}
\]
which combined with (\ref{V:BndsEst}) implies the lemma.
\end{proof}
 Now let
\begin{equation}\label{H:nonsmoothoperator}
G(v) := v - P_{[a,b]}\left(-\frac{1}{\alpha}p(y(v))\right), \mbox{ and
  } G_h(v) := v - P_{[a_h,b_h]}\left(-\frac{1}{\alpha}p_h(y_h(v))\right),
\end{equation}
where for given $v \in L^2(\Omega)$ the functions $p,p_h$ are defined through
\eqref{H:adjoint} and \eqref{H:adjointd}, respectively. It follows from the characterization of orthogonal projectors in real Hilbert spaces that the unique solutions
$u,u_h$ to \eqref{H:MP} and \eqref{H:MPd} are characterized by the equations
\begin{equation}\label{H:OE1}
G(u), \ G_h(u_h) = 0 \mbox{ in } {L^2(\Omega)}.
\end{equation}
These equations will be shown to be amenable to semi--smooth Newton methods as proposed in \cite{HIK03} and \cite{UM03}. We begin with formulating
\begin{Algorithm}\label{H:SSNM}{\rm (Semi--smooth Newton algorithm for \eqref{H:OE1})
    \begin{itemize}
    \item[] Start with $v \in {L^2(\Omega)}$ given. Do until convergence
    \item[] Choose $M\in\partial G_h(v)$.
    \item[] Solve $M\delta v = -G_h(v)$, $v := v+\delta v$.
    \end{itemize}
}
\end{Algorithm}
If we choose Jacobians $M\in\partial G_h(v)$ with $\|M^{-1}\|$ uniformly bounded throughout the iteration, and at the solution $u_h$ the function $G_h$ is $\partial G_h$-semismooth of order $\mu$, this algorithm is locally superconvergent of order $1+\mu$. 
 Although Algorithm \ref{H:SSNM} works on the infinite dimensional space $L^2(\Omega)$,
 it is possible to implement it numerically, as is shown subsequently.

\subsection{Semismoothness}
To apply the Newton algorithm, we need to confirm that the discretized operator $G_h$ is indeed semismooth. To establish this fact we rewrite $G_h$ in the form
\[
G_h(u)= u - \big(b-a\big)P_{[0,1]}\Big(\big(b-a\big)^{-1}\Big(-\frac{1}{\alpha}\big (S^*_h(S_h u-z)\big)-a\Big)\Big)+a
\]
and apply (\cite{UM03}, Theorem 5.2), with $P_{[0,1]}:\mathbb R\rightarrow\mathbb R$ taking the role of $\psi$. Here and in the following, for notational convenience we assume $a,b\in Y_h$, which is no restriction due to Lemma \ref{V:LemmaPertBnds}. The smoothing-operator $F:L^2\rightarrow L^q$ from \cite{UM03} in our case reads
\[
F(u) = \big(b-a\big)^{-1}\big(-\frac{1}{\alpha}\big (S^*_h(S_h u-z)\big)-a\big)\,.
\] 
We note that
\begin{itemize}
\item since we require $a,b\in L^\infty(\Omega)$, $b-a>\sigma>0$ with $\sigma\in \mathbb R$, both $(b-a)$ and $(b-a)^{-1}$ are in $L^\infty(\Omega)$ and the
pointwise multiplication by either $(b-a)$ or $(b-a)^{-1}$ is a continuous endomorphism in $L^p(\Omega)$ for any $p$.
\item the operator $F$ is differentiable with constant derivative for any $q\ge 1$. In fact, for sufficiently smooth domains $\Omega$,  the operators $S_h$ and $S_h^*$ map $L^2(\Omega)$ continuously into $H^2(\Omega)$, which is continuously embedded in $L^q(\Omega)$ for any $q\in [1,\infty]$.
\item $P_{[0,1]}:\mathbb R\rightarrow\mathbb R$ is $\partial P_{[0,1]}$-semismooth of order $1$, with 
\begin{equation}\label{V:genDiff}
\partial P_{[0,1]}(x)=\left\{\begin{array}{cl} 0&\textrm{if }x\notin\textrm{[0,1]}\\
1&\textrm{if }x\in\textrm{(0,1)}\\
\textrm{[0,1]}&\textrm{if }x=0\textrm{ or }x=1 \end{array}\right.\,.
\end{equation}
\item for piecewise linear elements the semismooth complementarity condition (5.3) in (\cite{UM03}, theorem 5.2) holds automatically with $\gamma=1$.
\end{itemize}

Thus we are in the position to apply (\cite{UM03},theorem 5.2) with $\alpha=1$ and $q_0>r=2$ and $\gamma=1$ and obtain
\begin{Theorem}
The function $G_h$ defined in (\ref{H:nonsmoothoperator}) is $\partial G_h$-semismooth of order $\mu<\frac{1}{3}$. There holds
\begin{equation*}
\partial G_h(v)w=w+\frac{1}{\alpha}\partial P_{[a,b]}\left(-\frac{1}{\alpha}\adj{v}\right)\cdot \big(S_h^*S_hw\big)\;,
\end{equation*}
where the application of the differential $\partial P_{[a,b]}$ and the multiplication by $S_h^*S_hw$ are pointwise operations a.e. in $\Omega$. 
\end{Theorem}
\begin{Remark}\label{R:MeshIndependence}
In \cite{HintermuellerUlbrich2004} the mesh independence of the superlinear convergence is stated.
Recent results from \cite{UM09} indicate semismoothness of $G$ of order $\frac{1}{2}$ as well as mesh independent q-superlinear convergence of the Newton algorithm of order $\frac{3}{2}$, if for example the modulus of the slope of $-\frac{1}{\alpha}p(y(\bar u))$ is bounded away from zero on the border of the active set, and if the mesh parameter $h$ is reduced appropriately. This is the key to our second globalization strategy proposed in Section \ref{sec:glob}
\end{Remark}

\subsection{Newton-Algorithm}

The generalized differential $\partial P_{[a,b]}$ can be defined analogously to (\ref{V:genDiff}) and the set-valued function $\partial P_{[a,b]}\big(-\frac{1}{\alpha}\adj{v}\big)$ contains the characteristic function $\chi_{\mathcal I(v)}$ of the inactive set
\[
\mathcal I(v)=\big\{\omega\in \Omega\,\big|\,\big(-\frac{1}{\alpha}\adj{v}\big)(\omega)\in\big(a(\omega),b(\omega)\big)\big\}\,.
\]
By $\chi^v$ we will denote synonymously the characteristic function $\chi_{\mathcal I(v)}$ as well as the self-adjoint endomorphism in $L^2(\Omega)$ given by the pointwise multiplication with $\chi_{\mathcal I(v)}$. 
The Newton-step in Algorithm \ref{H:SSNM} now takes the form
\begin{equation}\label{V:Newton_naiv}
\Big(I+\frac{1}{\alpha}\chi^v S_h^*S_h\Big)\delta v= - v + P_{[a,b]}\left(-\frac{1}{\alpha}\adj{v})\right)\,.
\end{equation}
To obtain an impression of the structure of the next iterate $v^+=v+\delta v$ we rewrite (\ref{V:Newton_naiv}) as
\[
v^+= P_{[a,b]}\left(-\frac{1}{\alpha}\adj{v}\right)-\frac{1}{\alpha}\chi^v S_h^*S_h\delta v\,.
\]
Since the range of $S^*_h$ is $Y_h$, the first addend is continuous and piecewise polynomial (of degree $k$) on a refinement $K_h$ of $T_h$.  The partition $K_h$ is obtained from $T_h$ by inserting nodes and edges along the boundary between the inactive set $\mathcal I(v)$ and the according active set, and in general contains simplices of higher order than $T_h$. The inserted edges are level sets of polynomials of order $\le k$ since we assume $a,b\in Y_h$.

 The second addend, involving the cut-off function $\chi^v$, is also piecewise polynomial of degree $k$ on $K_h$ but may jump along the edges not contained in $T_h$. 

Finally  $v^+$ lies in the following finite dimensional subspace of $L^2(\Omega)$
\[
Y_h^+=\left\{\chi^v\varphi_1+(1-\chi^v)\varphi_2\,\big|\;\varphi_1,\varphi_2\in Y_h\right\}= \text{span} \left(\{\phi_j\chi^v\}_{j=1}^n,\{\phi_j(1-\chi^v)\}_{j=1}^n\right)\,.
\]

The iterates generated by the Newton-algorithm can be represented exactly with about constant effort, since the number of inserted nodes varies only mildly from step to step, once the algorithm begins to converge. Furthermore the number of inserted nodes is bounded, see \cite{H03},\cite{H05}.

Since the Newton-increment $\delta v$ may have jumps along the borders of both the new and the old active and inactive sets, it is advantageous to compute $v^+$ directly, because $v^+$ lies in $Y_h^+$. To achieve an equation for $v^+$ we add $G_h'(v)v$ on both sides of (\ref{V:Newton_naiv}) to obtain
\begin{equation}\label{V:Newton}
\Big(I+\frac{1}{\alpha}\chi^v S_h^*S_h\Big)v^+= P_{[a,b]}\left(-\frac{1}{\alpha}\adj{v})\right)+\frac{1}{\alpha}\chi^v S_h^*S_h v\,,
\end{equation}
and reformulate Algorithm \ref{H:SSNM} as
\begin{Algorithm}[Newton Algorithm]\label{V:alg:Newton}
$\textrm{\hspace{4 cm}}$
    \begin{itemize}
    \item[] $v \in U$ given. Do until convergence
    \item[] Solve (\ref{V:Newton}) for $v^+$, $v := v^+$.
    \end{itemize}
\end{Algorithm}

\subsection{Computing the Newton-Step \ref{V:Newton}}

Since $v^+$ defined by (\ref{V:Newton}) is known on the active set $\mathcal A(v):=\Omega\setminus\mathcal I(V)$ it remains to compute $v^+$ on the inactive set. So we rewrite (\ref{V:Newton}) in terms of the unknown $\chi^vv^+$ by splitting $v^+$ as
\[
v^+=(1-\chi^v)v^++\chi^vv^+
\]
and obtain
\[
\Big(I+\frac{1}{\alpha}\chi^v S_h^*S_h\Big)\chi^vv^+= P_{[a,b]}\left(-\frac{1}{\alpha}\adj{v})\right)+\frac{1}{\alpha}\chi^v S_h^*S_h v
-\Big(I+\frac{1}{\alpha}\chi^v S_h^*S_h\Big)(1-\chi^v)v^+\,.
\]
As $(1-\chi^v)v^+$ is already known, we can restrict the latter equation to the inactive set $\mathcal I(v)$
\begin{equation}\label{V:RedSys}
\Big(\chi^v+\frac{1}{\alpha}\chi^v S_h^*S_h\chi^v\Big)v^+=\frac{1}{\alpha}\chi^vS_h^*z-\frac{1}{\alpha}\chi^v S_h^*S_h(1-\chi^v)v^+\,.
\end{equation}
On the left-hand side of (\ref{V:RedSys}) we have now a continuous, selfadjoint Operator on $L^2(\mathcal I^v)$, which is positive definite, because it is the restriction of the positive definite Operator $\left(I+\frac{1}{\alpha}\chi^vS_h^*S_h\chi^v\right)$ to $L^2(\mathcal I^v)$. 

Hence we are in the position to apply a CG-algorithm to solve (\ref{V:RedSys}). Moreover under the assumption of the first iterate lying in \[
Y_h^+\big|_{\mathcal I^v}=\left\{\chi^v\varphi\,\big|\;\varphi\in Y_h\right\}\,,
\]
as does the solution $\chi^vv^+$, the algorithm does not leave this space because of
\[
\left(I+\frac{1}{\alpha}\chi^vS_h^*S_h\chi^v\right)Y_h^+\big|_{\mathcal I^v}\subset Y_h^+\big|_{\mathcal I^v}
\]
and all CG-iterates lie in $Y_h^+\big|_{\mathcal I^v}$. These considerations lead to the following
\begin{Algorithm}[Solving (\ref{V:Newton})]\label{V:algSolve}
$\textrm{\hspace{3cm}}$
\begin{itemize}
  \item[] Compute the active and inactive sets $\mathcal A^v$ and $\mathcal I^v$.
  \item[] $\forall q\in \mathcal A^v$ set
  \[
  v^+(q)=P_{[a,b]}\left(-\frac{1}{\alpha}p_h(y_h(v))(q)\right)\,.
  \]
  \item[] Solve
  \[
  \Big(I+\frac{1}{\alpha}\chi^v S_h^*S_h\Big)\chi^vv^+=\frac{1}{\alpha}\chi^vS_h^*z-\frac{1}{\alpha}\chi^v S_h^*S_h(1-\chi^v)v^+
  \]
  for $\chi^vv^+$ by CG-iteration. By choosing a starting point in $Y_h^+\big|_{\mathcal I^v}$ one ensures that all iterates lie inside $Y_h^+\big|_{\mathcal I^v}$.
  \item[] $v^+=(1-\chi^v)v^++\chi^vv^+$.
\end{itemize}
\end{Algorithm}
We note that the use of this procedure in Algorithm \ref{V:alg:Newton} coincides with the active set strategy proposed in \cite{HIK03}.

\subsection{Globalization}\label{sec:glob}

Globalization of Algorithm \ref{V:alg:Newton} may require a damping step of the form
\[
v^+_\lambda = v+\lambda(v^+-v)
\]
with some $\lambda>0$. According to the considerations above, we have
\[
v^+_\lambda= (1-\lambda)v + \lambda \Big(P_{[a,b]}\left(-\frac{1}{\alpha}\adj{v}\right)-\frac{1}{\alpha}\chi^v S_h^*S_h\delta v\Big)\,.
\]
Unless $\lambda=1$ the effort of representing $v^+_\lambda$ will in general grow with every iteration of the algorithm, due to the jumps introduced in each step. 
This problem can be bypassed by focussing on the adjoint state $p_h(v)$ instead of the control $v$. In fact the function $\chi^v$ and thus also Equation (\ref{V:Newton}) do depend on $v$ only indirectly via the adjoint $p=p_h(v)=S_h^*(S_hv -z)$
\begin{equation}\label{V:GlobNewton}
\Big(I+\frac{1}{\alpha}\chi^p S_h^*S_h\Big)v^+= P_{[a,b]}\left(-\frac{1}{\alpha}p\right)+\frac{1}{\alpha}\chi^p( p+S_h^*z)\,.
\end{equation}
Now in each iteration the next full-step iterate $v^+$ is computed from (\ref{V:GlobNewton}). If damping is necessary, one computes $p^+_\lambda=p_h(v^+_\lambda)$ instead of $v_\lambda^+$.
In our (linear) setting the adjoint state $p^+_\lambda$ simply is a convex combination of $p=p_h(v)$ and $p^+=p_h(v^+)$
\[
p_\lambda ^+=\lambda p^+ +(1-\lambda)p\,,
\]
and unlike $v_\lambda^+$ the adjoint state $p_\lambda^+$ lies in the finite element space $Y_h$. Thus only a set of additional nodes according to the jumps of the most recent full-step iterate $v^+$ have to be managed, exactly as in the undamped case.
\begin{Algorithm}[Dampened Newton-Algorithm]\label{V:dampenedNewton}
$v \in U$ given.\\ Do until convergence
    \begin{itemize}
       \item[] Solve Equation (\ref{V:GlobNewton}) for $v^+$.
       \item[] Compute $p^+=p_h(y_h(v^+))$.
       \item[] Choose the damping-parameter $\lambda$. (for example by Armijo line search)
       \item[] Set $p:=p_\lambda ^+=\lambda p^+ +(1-\lambda)p$. 
   \end{itemize}
\end{Algorithm}
Algorithm \ref{V:alg:Newton} is identical to Algorithm \ref{V:dampenedNewton} without damping ($\lambda=1$).
\begin{Remark}
The above algorithm is equivalent to a dampened Newton algorithm applied to the equation
\begin{equation*}
p_h = S_h^*S_hP_{[a,b]}\left(-\frac{1}{\alpha}p_h\right)-S_h^*z\;,\qquad u := P_{[a,b]}\left(-\frac{1}{\alpha}p_h\right)\;.
\end{equation*}
\end{Remark}

Another approach, leading a globalization of Algorithm \ref{V:alg:Newton}, is to use some globalized, fully discrete scheme  and then perform \ref{V:alg:Newton} as a post processing step, compare also \cite{MeyerRoesch2004}.

Suppose $v_h$ is a discrete approximation to the optimal control $u$, such that $$\|v_h- u\|_{L^2(\Omega)}=\mathcal O(h)\,,$$ and let $ u_{hh}$ be its variationally discretized counterpart solving $(\mathbb P_{hh})$. Now, if the q-superlinear convergence of order $\frac{3}{2}$ of the Newton algorithm is mesh independent (see Remark \ref{R:MeshIndependence}),then  there exists a radius $\delta$ and a mesh parameter $h_0>0$, such that inside the ball $B_\delta(u_{hh})$ and for $h\le h_0$ the Newton iteration for $G_h$ converges q-superlinearly of order $\frac{3}{2}$ towards $u_{hh}$.

Let $\tilde u_{hh}$ be the second iterate of Algorithm \ref{V:alg:Newton} initialized with $v_h$. Then, for sufficiently small $h$, we have $v_h\in B_\delta(\bar u_h)$ and thus 
$$\|\tilde u_{hh} -u\|_{L^2(\Omega)}\le \|\tilde u_h -u_{hh}\|_{L^2(\Omega)}+ \|u_{hh} -u\|_{L^2(\Omega)} \le \|v_h -u_{hh}\|_{L^2(\Omega)}^{\frac{9}{4}}+\mathcal O(h^2)=\mathcal O(h^2)\,.$$

This motivates

\begin{Algorithm}[Post Processing]\label{V:PostProcessing}
$\textrm{\hspace{3cm}}$
    \begin{itemize}
    \item[] Solve the fully discretized optimization problem.
    \item[] Perform 2 steps of Algorithm \ref{V:algSolve}.
    \end{itemize}
\end{Algorithm}

\subsection{Global Convergence of the undamped Newton Algorithm}

It is not difficult to see, that the fixed-point equation for problem $(\mathbb P_{hh})$
\[
u_{hh}=P_{[a_h,b_h]}\left(-\frac{1}{\alpha}S_h^*(S_hu_{hh}-z)\right)
\]
can be solved by simple fixed-point iteration that converges globally for $\alpha>\|S_h\|_{L^2(\Omega),L^2(\Omega)}^2$, see \cite{H03,H05}. A similar global convergence result holds for the undamped Newton algorithm \ref{V:alg:Newton}
\begin{Lemma}\label{lemma:blob_conv}
The Newton algorithm \ref{V:alg:Newton} converges globally if $\alpha>\frac{4}{3}\|S\|^2$.
\end{Lemma}
\begin{proof}
See \cite{Dipl_Thesis_Vierling}.
\end{proof}

\section{Numerical examples}

We end this paper by illustrating our theoretical findings by
numerical examples. The first two examples are solved by Algorithm \ref{V:alg:Newton}, i.e. Algorithm \ref{V:dampenedNewton} without damping, making use of the global convergence property from Lemma \ref{lemma:blob_conv}. The third one involves a small parameter $\alpha=10^{-7}$ and is hence treated using the globalization strategy \ref{V:dampenedNewton} with Armijo line search. Finally the globalization \ref{V:PostProcessing} is applied at multiple parameters $\alpha$ and mesh parameters $h$.

As stopping criterion we require $\|P_{[a,b]}(-\frac{1}{\alpha}p_\lambda^+)-\bar u_h\|_{L^2(\Omega)}< 10^{-11}$ in Algorithm \ref{V:dampenedNewton}, using the a posteriori bound for admissible $v\in U_{ad}$
\[
\|v-\bar u_h\|_{L^2(\Omega)}\le \frac{1}{\alpha}\|\zeta\|_{L^2(\Omega)}\,,\qquad\zeta (\omega)=\left\{ \begin{array}{rl}
\,[\alpha v + p_h(v)]_{-} & \text{if }v(\omega)=a\\
\,[\alpha v + p_h(v)]_{+} & \text{if }v(\omega)=b\\
\alpha v + p_h(v)\phantom{]_{+}} & \text{if } a<v(\omega)<b
\end{array}\right.,
\]
presented in \cite{KrumbiegelRoesch2008} and \cite{TV09}.

\begin{Example}[Dirichlet]\label{V:Bsp:dirich}
We consider problem $(\mathbb P)$ in (\ref{H:MP}) with controls $u\in L^2(\Omega)$ on the unit square $\Omega=(0,1)^2$ with $a\equiv 0.3$ and $b\equiv 1$. Further we set
\[
z=-4\pi^2\alpha\sin(\pi x) \sin(\pi y)+(S\circ\imath) r\;,\,\textrm{ where } r=\min \big(1,\max\big(0.3,2\sin(\pi x) \sin(\pi y))\big)\big)\,.
\]
The choice of parameters implies a unique solution $\bar u=r$ to the continuous problem $(\mathbb{P})$.
\end{Example}

Throughout this section, solutions to the state equation are approximated by continuous, piecewise linear finite elements on a quasiuniform triangulation $T_h$ with maximal edge length $h>0$. The meshes are generated through regular refinement starting from the coarsest mesh. 
\begin{figure}[t]
%
%
\begin{psfrags}%
\psfragscanon%
%
\psfrag{s09}[b][b]{\color[rgb]{0,0,0}\setlength{\tabcolsep}{0pt}\begin{tabular}{c}Initial Guess\end{tabular}}%
\psfrag{s10}[rt][rt]{\color[rgb]{0,0,0}\setlength{\tabcolsep}{0pt}\begin{tabular}{r}x\end{tabular}}%
\psfrag{s11}[lt][lt]{\color[rgb]{0,0,0}\setlength{\tabcolsep}{0pt}\begin{tabular}{l}y\end{tabular}}%
\psfrag{s12}[b][b]{\color[rgb]{0,0,0}\setlength{\tabcolsep}{0pt}\begin{tabular}{c}u(x,y)\end{tabular}}%
\psfrag{s13}[b][b]{\color[rgb]{0,0,0}\setlength{\tabcolsep}{0pt}\begin{tabular}{c}Step 1\end{tabular}}%
\psfrag{s14}[rt][rt]{\color[rgb]{0,0,0}\setlength{\tabcolsep}{0pt}\begin{tabular}{r}x\end{tabular}}%
\psfrag{s15}[lt][lt]{\color[rgb]{0,0,0}\setlength{\tabcolsep}{0pt}\begin{tabular}{l}y\end{tabular}}%
\psfrag{s16}[b][b]{\color[rgb]{0,0,0}\setlength{\tabcolsep}{0pt}\begin{tabular}{c}u(x,y)\end{tabular}}%
\psfrag{s17}[b][b]{\color[rgb]{0,0,0}\setlength{\tabcolsep}{0pt}\begin{tabular}{c}Step 2\end{tabular}}%
\psfrag{s18}[rt][rt]{\color[rgb]{0,0,0}\setlength{\tabcolsep}{0pt}\begin{tabular}{r}x\end{tabular}}%
\psfrag{s19}[lt][lt]{\color[rgb]{0,0,0}\setlength{\tabcolsep}{0pt}\begin{tabular}{l}y\end{tabular}}%
\psfrag{s20}[b][b]{\color[rgb]{0,0,0}\setlength{\tabcolsep}{0pt}\begin{tabular}{c}u(x,y)\end{tabular}}%
\psfrag{s21}[b][b]{\color[rgb]{0,0,0}\setlength{\tabcolsep}{0pt}\begin{tabular}{c}Step 3\end{tabular}}%
\psfrag{s22}[rt][rt]{\color[rgb]{0,0,0}\setlength{\tabcolsep}{0pt}\begin{tabular}{r}x\end{tabular}}%
\psfrag{s23}[lt][lt]{\color[rgb]{0,0,0}\setlength{\tabcolsep}{0pt}\begin{tabular}{l}y\end{tabular}}%
\psfrag{s24}[b][b]{\color[rgb]{0,0,0}\setlength{\tabcolsep}{0pt}\begin{tabular}{c}u(x,y)\end{tabular}}%
%
\psfrag{x01}[t][t]{0}%
\psfrag{x02}[t][t]{0.1}%
\psfrag{x03}[t][t]{0.2}%
\psfrag{x04}[t][t]{0.3}%
\psfrag{x05}[t][t]{0.4}%
\psfrag{x06}[t][t]{0.5}%
\psfrag{x07}[t][t]{0.6}%
\psfrag{x08}[t][t]{0.7}%
\psfrag{x09}[t][t]{0.8}%
\psfrag{x10}[t][t]{0.9}%
\psfrag{x11}[t][t]{1}%
\psfrag{x12}[t][t]{0}%
\psfrag{x13}[t][t]{0.5}%
\psfrag{x14}[t][t]{1}%
\psfrag{x15}[t][t]{0}%
\psfrag{x16}[t][t]{0.5}%
\psfrag{x17}[t][t]{1}%
\psfrag{x18}[t][t]{0}%
\psfrag{x19}[t][t]{0.5}%
\psfrag{x20}[t][t]{1}%
\psfrag{x21}[t][t]{0}%
\psfrag{x22}[t][t]{0.1}%
\psfrag{x23}[t][t]{0.2}%
\psfrag{x24}[t][t]{0.3}%
\psfrag{x25}[t][t]{0.4}%
\psfrag{x26}[t][t]{0.5}%
\psfrag{x27}[t][t]{0.6}%
\psfrag{x28}[t][t]{0.7}%
\psfrag{x29}[t][t]{0.8}%
\psfrag{x30}[t][t]{0.9}%
\psfrag{x31}[t][t]{1}%
\psfrag{x32}[t][t]{0}%
\psfrag{x33}[t][t]{0.5}%
\psfrag{x34}[t][t]{1}%
%
\psfrag{v01}[r][r]{0}%
\psfrag{v02}[r][r]{0.2}%
\psfrag{v03}[r][r]{0.4}%
\psfrag{v04}[r][r]{0.6}%
\psfrag{v05}[r][r]{0.8}%
\psfrag{v06}[r][r]{1}%
\psfrag{v07}[r][r]{0}%
\psfrag{v08}[r][r]{0.5}%
\psfrag{v09}[r][r]{1}%
\psfrag{v10}[r][r]{0}%
\psfrag{v11}[r][r]{0.5}%
\psfrag{v12}[r][r]{1}%
\psfrag{v13}[r][r]{0}%
\psfrag{v14}[r][r]{0.5}%
\psfrag{v15}[r][r]{1}%
\psfrag{v16}[r][r]{0}%
\psfrag{v17}[r][r]{0.2}%
\psfrag{v18}[r][r]{0.4}%
\psfrag{v19}[r][r]{0.6}%
\psfrag{v20}[r][r]{0.8}%
\psfrag{v21}[r][r]{1}%
\psfrag{v22}[r][r]{0}%
\psfrag{v23}[r][r]{0.5}%
\psfrag{v24}[r][r]{1}%
%
\psfrag{z01}[r][r]{0.3}%
\psfrag{z02}[r][r]{1}%
\psfrag{z03}[r][r]{1.5}%
\psfrag{z04}[r][r]{0.3}%
\psfrag{z05}[r][r]{1}%
\psfrag{z06}[r][r]{1.5}%
\psfrag{z07}[r][r]{0}%
\psfrag{z08}[r][r]{0.5}%
\psfrag{z09}[r][r]{1}%
\psfrag{z10}[r][r]{0.3}%
\psfrag{z11}[r][r]{1}%
\psfrag{z12}[r][r]{1.5}%
%
\resizebox{16cm}{!}{\includegraphics{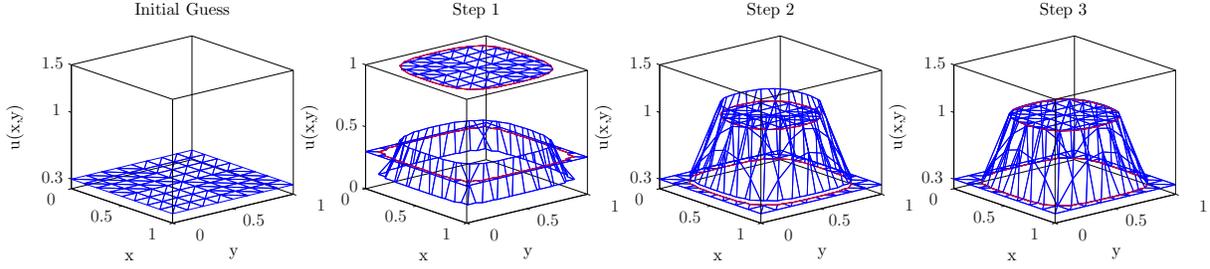}}%
\end{psfrags}%
%
\caption{The first four Newton-iterates for Example \ref{V:Bsp:dirich} (Dirichlet) with parameter $\alpha=0.001$}\label{V:fig:dirich}
\end{figure}
\begin{table}[b]
\center
\begin{tabular}{l|c|c|c|c|c|c}
mesh param. $h$ & $ERR$ & $ERR_{\infty}$ & $EOC$ & $EOC_{\infty}$ & Iterations & Quality\\
\hline
$\sqrt{2}/16$    & 2.5865e-03 & 1.2370e-02 & 1.95 & 1.79 & 4 & 2.16e-15\\
$\sqrt{2}/32$    & 6.5043e-04 & 3.2484e-03 & 1.99 & 1.93 & 4 & 2.08e-15\\
$\sqrt{2}/64$    & 1.6090e-04 & 8.1167e-04 & 2.02 & 2.00 & 4 & 2.03e-15\\
$\sqrt{2}/128$  & 4.0844e-05 & 2.1056e-04 & 1.98 & 1.95 & 4 & 1.99e-15\\
$\sqrt{2}/256$  & 1.0025e-05 & 5.3806e-05 & 2.03 & 1.97 & 4 & 1.69e-15\\
$\sqrt{2}/512$  & 2.5318e-06 & 1.3486e-05 & 1.99 & 2.00 & 4 & 1.95e-15
\end{tabular}
\caption{$L^2$- and $L^\infty$-error development for Example \ref{V:Bsp:dirich} (Dirichlet)}\label{V:tab:dirich}
\end{table}

As discussed in Section \ref{H:sec:FEsemiD}, problem $(\mathbb{P}_{hh})$ admits a unique solution $\bar u_h$ and we have
\[
\|\bar u_h-\bar u\|_{L^2(\Omega)}=O(h^2)
\]
as $h\rightarrow 0$. There also holds nearly quadratic convergence in $L^\infty(\Omega)$ 
\[
\|\bar u_h-\bar u\|_{L^\infty(\Omega)}=O(|\log(h)|^{\frac{1}{2}}h^2)
\]
for domains $\Omega\subset\mathbb R^2$, see \cite{H03}.
Both convergence rates are observed in Table \ref{V:tab:dirich}, that shows the $L^2$- and the $L^\infty$-errors together with the corresponding experimental orders of convergence
\[
EOC_i = \frac{\ln ERR(h_{i-1})-\ln ERR(h_i)}{\ln(h_{i-1})-\ln(h_i)}
\]
for Example \ref{V:Bsp:dirich}. Lemma \ref{lemma:blob_conv} ensures global convergence of the undamped Algorithm \ref{V:alg:Newton} only for $\alpha>1/(3\pi^4)\simeq 0.0034$, but it is still observed for $\alpha=0.001$. 

The algorithm is initialized with $v_0\equiv 0.3$. The resulting number of Newton steps as well as the value of $\zeta/\alpha$ for the computed solution are also given in Table \ref{V:tab:dirich}.

Figure \ref{V:fig:dirich} shows the Newton iterates, active and inactive sets are very well distinguishable, the jumps along their frontier can be observed.

\begin{table}[t]
\center
\begin{tabular}{l|c|c|c|c|c|c}
mesh param. $h$ & $ERR$ & $ERR_{\infty}$ & $EOC$ & $EOC_{\infty}$ & Iterations & Quality\\
\hline
$\sqrt{2}/16$    & 3.9866e-03 & 1.1218e-02 & 1.94 & 1.74 & 3 & 1.81e-12\\
$\sqrt{2}/32$    & 1.0025e-03 & 3.2332e-03 & 1.99 & 1.79 & 3 & 2.31e-12\\
$\sqrt{2}/64$    & 2.5188e-04 & 8.4398e-04 & 1.99 & 1.94 & 3 & 9.74e-13\\
$\sqrt{2}/128$  & 6.2936e-05 & 2.1856e-04 & 2.00 & 1.95 & 3 & 9.37e-13\\
$\sqrt{2}/256$  & 1.5740e-05 & 5.5223e-05 & 2.00 & 1.99 & 3 & 8.91e-13\\
$\sqrt{2}/512$  & 3.9346e-6 & 1.3928e-05 & 2.00 & 2.00 & 3 & 8.86e-13
\end{tabular}
\caption{Development of the error in Example \ref{V:Bsp:Neumann} (Neumann)}\label{tab:neum}
\end{table}

Next we demonstrate another Example, out theory may also be applied to.
\begin{Example}[Neumann]\label{V:Bsp:Neumann}
We next consider an elliptic problem with Neumann boundary conditions
\begin{equation*}
\begin{split}
 -\Delta y + y&= u \qquad \textrm{in $\Omega$}\;, \\
  \partial_n y &= 0 \qquad \textrm{on $\partial\Omega$}\;,
\end{split}
\end{equation*}
on $\Omega=(0,1)^2$, with a similar discrete setting as in the previous example. It then is clear, how $(\mathbb P)$ and $(\mathbb P_{hh})$ have to be understood.
We set $\alpha=1$ and choose
\[
z=-2(2\pi^2+1)\alpha\cos(\pi x) \cos(\pi y)+(S\circ\imath)r\;,\,\textrm{ with } r=\min \big(1,\max\big(-1,2\cos(\pi x) \cos(\pi y)\big)\big)
\]
and bounds $a\equiv-1$ and $b\equiv 1$. The optimal control to the continuous problem is $\bar u = r$. 
\end{Example}
\begin{figure}[b]
%
%
\begin{psfrags}%
\psfragscanon%
%
\psfrag{s09}[b][b]{\color[rgb]{0,0,0}\setlength{\tabcolsep}{0pt}\begin{tabular}{c}Initial Guess\end{tabular}}%
\psfrag{s10}[rt][rt]{\color[rgb]{0,0,0}\setlength{\tabcolsep}{0pt}\begin{tabular}{r}x\end{tabular}}%
\psfrag{s11}[lt][lt]{\color[rgb]{0,0,0}\setlength{\tabcolsep}{0pt}\begin{tabular}{l}y\end{tabular}}%
\psfrag{s12}[b][b]{\color[rgb]{0,0,0}\setlength{\tabcolsep}{0pt}\begin{tabular}{c}u(x,y)\end{tabular}}%
\psfrag{s13}[b][b]{\color[rgb]{0,0,0}\setlength{\tabcolsep}{0pt}\begin{tabular}{c}Step 1\end{tabular}}%
\psfrag{s14}[rt][rt]{\color[rgb]{0,0,0}\setlength{\tabcolsep}{0pt}\begin{tabular}{r}x\end{tabular}}%
\psfrag{s15}[lt][lt]{\color[rgb]{0,0,0}\setlength{\tabcolsep}{0pt}\begin{tabular}{l}y\end{tabular}}%
\psfrag{s16}[b][b]{\color[rgb]{0,0,0}\setlength{\tabcolsep}{0pt}\begin{tabular}{c}u(x,y)\end{tabular}}%
\psfrag{s17}[b][b]{\color[rgb]{0,0,0}\setlength{\tabcolsep}{0pt}\begin{tabular}{c}Step 2\end{tabular}}%
\psfrag{s18}[rt][rt]{\color[rgb]{0,0,0}\setlength{\tabcolsep}{0pt}\begin{tabular}{r}x\end{tabular}}%
\psfrag{s19}[lt][lt]{\color[rgb]{0,0,0}\setlength{\tabcolsep}{0pt}\begin{tabular}{l}y\end{tabular}}%
\psfrag{s20}[b][b]{\color[rgb]{0,0,0}\setlength{\tabcolsep}{0pt}\begin{tabular}{c}u(x,y)\end{tabular}}%
\psfrag{s21}[b][b]{\color[rgb]{0,0,0}\setlength{\tabcolsep}{0pt}\begin{tabular}{c}Step 3\end{tabular}}%
\psfrag{s22}[rt][rt]{\color[rgb]{0,0,0}\setlength{\tabcolsep}{0pt}\begin{tabular}{r}x\end{tabular}}%
\psfrag{s23}[lt][lt]{\color[rgb]{0,0,0}\setlength{\tabcolsep}{0pt}\begin{tabular}{l}y\end{tabular}}%
\psfrag{s24}[b][b]{\color[rgb]{0,0,0}\setlength{\tabcolsep}{0pt}\begin{tabular}{c}u(x,y)\end{tabular}}%
%
\psfrag{x01}[t][t]{0}%
\psfrag{x02}[t][t]{0.1}%
\psfrag{x03}[t][t]{0.2}%
\psfrag{x04}[t][t]{0.3}%
\psfrag{x05}[t][t]{0.4}%
\psfrag{x06}[t][t]{0.5}%
\psfrag{x07}[t][t]{0.6}%
\psfrag{x08}[t][t]{0.7}%
\psfrag{x09}[t][t]{0.8}%
\psfrag{x10}[t][t]{0.9}%
\psfrag{x11}[t][t]{1}%
\psfrag{x12}[t][t]{0}%
\psfrag{x13}[t][t]{0.5}%
\psfrag{x14}[t][t]{1}%
\psfrag{x15}[t][t]{0}%
\psfrag{x16}[t][t]{0.5}%
\psfrag{x17}[t][t]{1}%
\psfrag{x18}[t][t]{0}%
\psfrag{x19}[t][t]{0.5}%
\psfrag{x20}[t][t]{1}%
\psfrag{x21}[t][t]{0}%
\psfrag{x22}[t][t]{0.1}%
\psfrag{x23}[t][t]{0.2}%
\psfrag{x24}[t][t]{0.3}%
\psfrag{x25}[t][t]{0.4}%
\psfrag{x26}[t][t]{0.5}%
\psfrag{x27}[t][t]{0.6}%
\psfrag{x28}[t][t]{0.7}%
\psfrag{x29}[t][t]{0.8}%
\psfrag{x30}[t][t]{0.9}%
\psfrag{x31}[t][t]{1}%
\psfrag{x32}[t][t]{0}%
\psfrag{x33}[t][t]{0.5}%
\psfrag{x34}[t][t]{1}%
%
\psfrag{v01}[r][r]{0}%
\psfrag{v02}[r][r]{0.2}%
\psfrag{v03}[r][r]{0.4}%
\psfrag{v04}[r][r]{0.6}%
\psfrag{v05}[r][r]{0.8}%
\psfrag{v06}[r][r]{1}%
\psfrag{v07}[r][r]{0}%
\psfrag{v08}[r][r]{0.5}%
\psfrag{v09}[r][r]{1}%
\psfrag{v10}[r][r]{0}%
\psfrag{v11}[r][r]{0.5}%
\psfrag{v12}[r][r]{1}%
\psfrag{v13}[r][r]{0}%
\psfrag{v14}[r][r]{0.5}%
\psfrag{v15}[r][r]{1}%
\psfrag{v16}[r][r]{0}%
\psfrag{v17}[r][r]{0.2}%
\psfrag{v18}[r][r]{0.4}%
\psfrag{v19}[r][r]{0.6}%
\psfrag{v20}[r][r]{0.8}%
\psfrag{v21}[r][r]{1}%
\psfrag{v22}[r][r]{0}%
\psfrag{v23}[r][r]{0.5}%
\psfrag{v24}[r][r]{1}%
%
\psfrag{z01}[r][r]{-1}%
\psfrag{z02}[r][r]{0}%
\psfrag{z03}[r][r]{1}%
\psfrag{z04}[r][r]{0}%
\psfrag{z05}[r][r]{0.5}%
\psfrag{z06}[r][r]{1}%
\psfrag{z07}[r][r]{1.5}%
\psfrag{z08}[r][r]{-1}%
\psfrag{z09}[r][r]{0}%
\psfrag{z10}[r][r]{1}%
\psfrag{z11}[r][r]{0}%
\psfrag{z12}[r][r]{0.5}%
\psfrag{z13}[r][r]{1}%
\psfrag{z14}[r][r]{1.5}%
\psfrag{z15}[r][r]{-3}%
\psfrag{z16}[r][r]{-2}%
\psfrag{z17}[r][r]{-1}%
\psfrag{z18}[r][r]{0}%
\psfrag{z19}[r][r]{1}%
\psfrag{z20}[r][r]{-1}%
\psfrag{z21}[r][r]{0}%
\psfrag{z22}[r][r]{1}%
\psfrag{z23}[r][r]{0}%
\psfrag{z24}[r][r]{0.5}%
\psfrag{z25}[r][r]{1}%
\psfrag{z26}[r][r]{1.5}%
%
\resizebox{16cm}{!}{\includegraphics{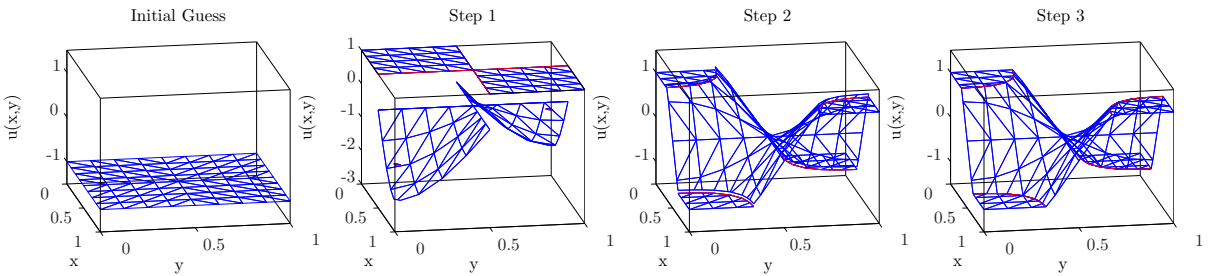}}%
\end{psfrags}%
%

\caption{The first steps of the Newton-algorithm for Example \ref{V:Bsp:Neumann} (Neumann) with $\alpha=1$.}\label{V:fig:Neum}
\end{figure}
For $\alpha =1 $ the undamped iteration still converges globally, although the solution operator has norm $\|S\|=1$ as an endomorphism in $L^2(\Omega)$. The predicted convergence properties  and the stopping criterion are the same as above; Algorithm \ref{V:dampenedNewton} is initialized by $v_0\equiv -1$. The first four steps of the iteration are displayed in Figure \ref{V:fig:Neum} and the behaviour of the approximation error between the exact and the semidiscrete solution, as well as the number of iterations and the final value of $\zeta/\alpha$, is shown in Table \ref{tab:neum}.

The Algorithm has also been implemented successfully for parabolic discontinuous Galerkin discretized problems as well as elliptic problems with Lavrentiev-regularized state constraints.\\[0.5ex]

To demonstrate Algorithm \ref{V:dampenedNewton} with damping we again consider Example \ref{V:Bsp:dirich}, this time with $\alpha=10^{-7}$. We choose 
$$MF(p) =   \left\|p-S_h^*S_hP_{[a,b]}\left(-\frac{1}{\alpha}p\right)+S_h^*z\right\|^2_{L^2(\Omega)}\;,$$ 
as merit function governing the step size of the algorithm. Again we use the same stopping criterion as in the previous examples.

Table \ref{tab:glob} shows errors and the number of iterations for different mesh parameters $h$ at a smoothing parameter $\alpha=10^{-7}$. To compare the number of iterations we choose a common initial guess $u_0 \equiv 1$. The number of iterations appears to be independent of $h$.

Finally, to demonstrate the efficiency of Algorithm \ref{V:PostProcessing}, the EOC in the $L^2(\Omega)$-norm is plotted in table \ref{tab:PostProc}. The disturbances that can be observed for smaller parameter $\alpha$ indicate the decay of the environment of q-superlinear convergence with decreasing $\alpha$.

\begin{table}[t]
\center
\begin{tabular}{l|c|c|c|c|c}
mesh param. $h$ & $ERR$ & $ERR_{\infty}$ & $EOC$ & $EOC_{\infty}$ & Iterations\\
\hline
$\sqrt{2}/2$      & 1.1230e-01 & 3.0654e-01 & -        & -       & 11\\
$\sqrt{2}/4$      & 3.8398e-02 & 1.4857e-01 & 1.55 & 1.04 & 22\\
$\sqrt{2}/8$      & 9.8000e-03 & 4.4963e-02 & 1.97 & 1.72 & 19\\
$\sqrt{2}/16$    & 1.7134e-03 & 1.2316e-02 & 2.52 & 1.87 & 20\\
$\sqrt{2}/32$    & 4.0973e-04 & 2.8473e-03 & 2.06 & 2.11 & 33\\
$\sqrt{2}/64$    & 8.2719e-05 & 6.2580e-04 & 2.31 & 2.19 & 17\\
$\sqrt{2}/128$  & 2.0605e-05 & 1.4410e-04 & 2.01 & 2.12 & 20\\
$\sqrt{2}/256$  & 4.7280e-06 & 4.6075e-05 & 2.12 & 1.65 & 19\\
$\sqrt{2}/512$  & 1.1720e-06 & 1.0363e-05 & 2.01 & 2.15 & 18
\end{tabular}
\caption{Development of the error in Example \ref{V:Bsp:dirich} (Dirichlet) for $\alpha=10^{-7}$.}\label{tab:glob}
\end{table}

\begin{table}[b]
\center
\begin{tabular}{|l|c|c|c|c|c|c|c|c|c|}
\hline
\backslashbox[2cm]{\phantom{M}$h$}{\\$\alpha$\phantom{M}} & $1$ & $10^{-1}$ & $10^{-2}$ & $10^{-3}$ & $10^{-4}$&$10^{-5}$&$10^{-6}$&$10^{-7}$&$10^{-8}$\\
\hline
$\sqrt{2}/4$     &     1.36    & 1.37   & 1.42   & 1.81&          1.81  &         1.56   &       1.24  &        0.96 &-4.07\\
$\sqrt{2}/8$     &     1.78    & 1.78   & 1.75   & 1.53&          1.94  &        2.41    &      2.41   &      -6.49 &-2.67\\
$\sqrt{2}/16$   &     1.95    & 1.96   & 1.97   & 2.09&           2.35 &         2.68   &       2.01  &         9.41& 2.96\\
$\sqrt{2}/32$   &     2.01    &  2.01  & 2.00   & 1.94&          1.83  &       1.78     &     2.95    &      4.40 &10.43\\
$\sqrt{2}/64$   &     1.99    &  1.99  & 1.99   & 2.02&          2.09  &        2.19    &      2.13   &      1.84  &-1.43\\
$\sqrt{2}/128$ &     2.00    &  2.00  & 2.00   & 1.99&          1.97  &        1.92    &      2.01   &       2.50 & 7.05\\
$\sqrt{2}/256$ &     2.00    &   2.00 & 2.00   & 2.01&          2.04  &        2.08    &      2.11   &      2.19 &2.34\\
\hline
\end{tabular}
\caption{EOC of Algorithm \ref{V:PostProcessing} applied to Example \ref{V:Bsp:dirich} (Dirichlet) for different $\alpha$.}\label{tab:PostProc}
\end{table}


\section*{Acknowledgements} 
The first author gratefully acknowledges the support of the DFG
Priority Program 1253 entitled Optimization With Partial
Differential Equations. We also thank Andreas G\"unther for some fruitful discussions.

\end{document}